\documentclass[12pt]{amsart}
\usepackage{amsfonts}

\theoremstyle{definition}

\theoremstyle{remark}

\begin{document}

\begin{center}

 {\bf MULTIDIMENSIONAL DILATION OPERATORS, BOYD AND SHIMOGAKI INDICES OF
BILATERAL WEIGHT GRAND LEBESQUE SPACES} \par

\vspace{3mm}

{\bf E. Ostrovsky.} \ Email: \ galo@list.ru \par

\vspace{2mm}

{\bf L. Sirota.} \ Email: \ sirota@zahav.net.il \par

\vspace{3mm}

 Abstract \par

\vspace{2mm}

 In this paper we compute the norm of dilation operators,
multidimensional Boyd's and Shimogaki's indices in the Bilateral
 Grand Lebesgue Spaces and consider some applications. \\

\end{center}

\vspace{3mm}


\vspace{3mm}


\vspace{3mm}

\section{Introduction. Statement of problem}

\vspace{3mm}

 Let $ X $ be a $ d \ $ dimensional positive Euclidean subspace with a "multiplicative" decomposition of a view:

$$
X = R_+^{d(1)} \times R_+^{d(2)} \times \ldots \times R_+^{d(k)}, k = 1,2,\ldots,
$$

so that $ d = d(1) + d(2) + \ldots + d(k). $ We will write for
$ \vec{x} = x \in X $

$$
\vec{x} = x = (x(d(1)), x(d(2)), \ldots, x(d(k)), \ x(d(j)) = \vec{x}(d(j))
\in R_+^{d(j)}, \ j = 1,2, \ldots, k.
$$

 The space $ X $ is equipped by the usually Borelian sigma \ - \ algebra
$ \Sigma $ and by the non \ - \ negative weight measure

$$
\mu(V) = \mu_W(V) = \int_V W(\vec{x}) \ d \vec{x},  \ V \in \Sigma,
$$

$$
W(\vec{x}) = \prod_{r=1}^k W_r( \vec{x}( d(r))) ,
$$
where all the non \ - \ trivial functions $ W_r(\cdot) $ are continuous in
the domain $ X \setminus \{ 0 \}, $ non \ - \ negative and
homogeneous of an order $ \theta(r): \ \forall \lambda > 0 \ \Rightarrow $

$$
W_r(\lambda \ x(d(r))) = \lambda^{\theta(r)} W_r(x(d(r))).
$$
Notation:
$$
\vec{\theta} = \theta = \{\theta(1), \theta(2), \ldots, \theta(k) \}.
$$

 For $a$ and $b$ constants, $1 \le a < b \le \infty,$ let $\psi =
\psi(p) = \psi(p; a,b),$ $ p \in (a,b),$ be a continuous on the {\it open} interval $ (a,b) $ positive: $ \psi(p) \ge 1 $
function such that $\psi(a + 0)$ and $\psi(b-0)$ exist, both may be infinite,
and postulate that $ \psi(b-0)= \infty. $ \par

The class such a functions we will denote $ E\Psi = E\Psi(a,b);
 \ E\Psi = \{\psi \}. $\par

 The Bilateral Grand Lebesgue space (in notation BGL, BGLS)
$  G_X(\mu; \psi; a,b) =
G_X(\psi;a,b) = G(\psi; a,b) = G(\psi) $ is the space of all the measurable
functions $ h: X \to{\mathbb R} $ endowed with the norm

\begin{eqnarray*}
||h||G(\psi) \stackrel{def}{=}\sup_{p \in (a,b)}|h|_p/\psi(p),
 \quad |h|_p = |h|_{p, \mu} = \left[\int_X |h(x)|^p \ d\mu(x) \right]^{1/p}.
\end{eqnarray*}

 By definition, $ h(\cdot) \in G_X(\mu; \psi, a,b) = G(\psi) $ if and only
if  $ ||h||G(\psi) < \infty. $ \par

\vspace{3mm}

{\bf Proposition 1.}  {\it The BGL spaces are Banach functional spaces; moreover, they are rearrangement invariant (r.i.) spaces. }\par

\vspace{3mm}

{\bf Proof.} We must prove only that the $ G(\psi) $ space satisfies the Fatou property, since all the other properties are evident. \par
 Recall at first that the Fatou property of some r.i. space $ G $ over
source triplet $ (X, \Sigma, \mu ) $ denotes that for arbitrary
non-increasing sequence of non-negative functions $ \{  f_n \}
= \{ f_n(x), \ x \in X  \} $ belonging to the space $ G $ and such that  as
$ n \uparrow \infty $

$$
f_n(x) \uparrow f(x), \ \sup_n ||f_n||G  < \infty
$$
it follows

$$
||f_n||G \uparrow ||f||G.
$$
 Let $ G = G(\psi) $ and suppose the sequence of the measurable
functions $ \{f_n \} = \{f_n: X \to R \} $ satisfies our condition.
As long as the space $ L_p(X, \mu) $ satisfies the Fatou property, we have:

$$
\sup_n ||f_n||G(\psi) = \sup_n \sup_{p \in (a,b)} [ |f_n|_p/\psi(p) ] =
$$

$$
\sup_{p \in (a,b)} \sup_n [|f_n|_p/\psi(p) ] = \sup_{p \in (a,b)}[|f|_p
/\psi(p)] = ||f||G(\psi),
$$
Q.E.D. \par
 Notice that in this proof we do not use the concrete view of the measure
$ \mu; $ we assume only that the measure $ \mu $ is sigma \ - \ finite. \par

 The  $G(\psi)$ spaces with $ \mu(X) = 1$ appeared in \cite{KO11};
it was proved that in this case each $ G(\psi) $ space coincides
with certain exponential Orlicz space, up to norm equivalence.
Partial cases of these spaces were intensively studied, in
particular, their associate spaces, fundamental functions
$\phi(G(\psi; a,b);\delta),$ Fourier and singular operators,
conditions for convergence and compactness, reflexivity and
separability, martingales in these spaces, etc.; see, e.g.,
\cite{BS1},\cite{CM3}-\cite{JM9},\cite{KM10, Holon,OS14}. These
spaces are also Banach and moreover rearrangement invariant
(r.i.).\par
 Some classical properties of these spaces  (Sobolev embedding theorems,
convolution operators etc.) was partially investigated in \cite{LOS}. \par

The BGLS norm estimates, in particular, Orlicz norm estimates for
measurable functions, e.g., for random variables are used in PDE
\cite{Fio4}-\cite{IKO8}, probability in Banach spaces
\cite{LedTal18}, in the modern non-parametrical statistics, for
example, in the so-called regression problem \cite{Os13,OS14}.

We are going to deal not with all functions $\psi$ described above
but with an essential subset of these functions satisfying certain
natural conditions.

Let again $ a \ge 1,$ $ b\in (a,\infty],$ and let $ \psi = \psi(p)
$ be a positive continuous function on the {\it open} interval $
(a,b) $ such that  there exists a measurable function $f: X
\to{\mathbb R} $ for which

\begin{eqnarray}\label{no0}
f(\cdot) \in \cap_{p \in (a,b) } L_p(X,\mu) = \cap_{p \in (a,b) }L_p,
\ \psi(p) = |f|_p, \ p \in (a,b).
\end{eqnarray}

and such that $\max\{\psi(a+0), \psi(b-0)\} = \infty $ and in the case $ b = \infty $ we define $ \psi(b-0) = \lim_{p \to \infty} \psi(p) $
and postulate again that $ \psi(b-0)= \infty. $ \par

  We say that the equality (\ref{no0}) and the function $  f(\cdot) $ from
(\ref{no0}) is the {\it representation} of the function $ \psi. $
The existence of representation implies, by the way, the
log-convexity of $\psi.$

We denote the subset of all the functions $\psi$ having
representation by $ \Psi= \Psi(a,b).$ For complete description of
these functions see, for example, (\cite[p.p. 21-27]{Os13},
\cite{OS14}). \par

 Another definitions.
We will say as usually  ( \cite{BS1}, p. 14 \ - \ 16)
that the function $ f \in G(\psi), \ \psi \in E\Psi $ has
{\it absolutely continuous norm and write }$ f \in GA(\psi), $
{\it if}
$$
 \lim_{\delta \to 0} \sup_{A: \mu(A) \le \delta} ||f \ I_A||G(\psi)=0.
$$

 We denote by $  G^o = G^o(\psi) = G^o_X(\psi), \ \psi \in E \Psi $ the closed subspace of $ G(\psi), $ consisting on all the functions $ f, $ satisfying the following condition:
$$
\lim_{\psi(p) \to \infty} |f|_p/\psi(p) = 0;
$$
 and  denote by $ GB =  GB(\psi) $ the closed span in the norm $ G(\psi) $
the set of all bounded: $ vraimax_x |f(x)| < \infty $
 measurable functions with finite support: $ \mu(supp \ |f|) < \infty. $ \par

The subspaces $ GA(\psi), GB(\psi), G^0(\psi) $ are closed r.i.
subspaces of the space $ G(\psi). $ \par

It is proved in \cite{OS14}, \cite{Holon} that  if $ \psi \in \Psi $ then
$$
G(\psi) \ne G^0(\psi) = GB(\psi) = GA(\psi) \subset G(\psi).
$$

 We note also that the $ G(\psi) $ spaces are also interpolation spaces
(the so-called  $\Sigma$-spaces), see \cite{BS1},\cite{CM3}-\cite
{JM9},\cite{Holon}, \cite{OS14}, \cite{Nik} etc. However, we hope that our direct representation of these spaces is of certain convenience in both
theory and applications. A natural question arises what happens if
the spaces other than $L_p$ are used in the definition. Indeed,
this is possible and might be of interest, but, for example, using
Lorenz spaces in this capacity leads to the same object (see
\cite{OS14} ).\par

 We will use widely further the notion of {\it fundamental function}
$ \phi(G, \delta), \ \delta \in (0,\infty) $ of the arbitrary r.i. space
$ G $ over the triplet $ (X, \Sigma, \mu). $ Recall that  by definition

$$
 \phi(G, \delta) = || \ I(A) \ ||G, \ \mu(A) = \delta, \ \delta \in
[0, \infty)
$$
and $ I(A) = I(A,x) = 1, x \in A, \ I(A) = I(A,x) = 0, x \notin A. $ \par
 This notion play a very important role in the theory of interpolation of
operators, theory of Fourier series, theory of approximation etc. See, for example, \cite{LOS}, \cite{BS1}, \cite{Nik} etc.\par

 {\bf Lemma 1. }
$$
 \psi(\cdot) \in E \Psi \ \Rightarrow \ \phi(G(\psi), 0+ ) = 0.
$$
{\bf Proof.}  Assume at first that $ b < \infty; $ then

$$
\phi(G(\psi), \delta) \le \delta^{1/b} \to 0, \ \delta \to 0.
$$

 Therefore, it remains only to consider the case $ b = \infty. $ \par
 Recall that in this case we suppose $ \psi(\infty) = \infty. $ \par
 Let $ \epsilon = const \in (0, 1) $ be arbitrary "small" number. There
exists  a number $ Q > 1 $ for which

$$
\forall p > Q \ \Rightarrow \psi(p) > 1/\epsilon.
$$
 We have for all sufficient small values $ \delta: \ \delta \in (0, \delta_0),
\delta_0 = \delta_0(\epsilon) \in (0,1):$
$$
\phi(G(\psi), \delta) \le \delta^{1/Q} + \epsilon < 2 \epsilon.
$$
As a consequence: \par
{\bf Lemma 2. } Assume that the metric space $ \Sigma $ under the distance

$$
\rho(A,B) = \arctan( \mu ( (A \setminus B) \ \cup (B \setminus A)) ), \ A,B \in
\Sigma
$$
is separable. Then the r.i. space $ G^o(\psi) (= GA(\psi) = GB(\psi)) $ is
also separable. \par

{\bf Lemma 3.} Let  $ \psi(\cdot) \in \Psi(a,b). $ Then

$$
\lim_{s \to 0 +} \frac{ \log \phi (G(\psi),s) }{ \log s } = \frac{1}{b};
$$

$$
\lim_{s \to \infty} \frac{ \log \phi ( G (\psi),s ))}{ \log s } = \frac{1}{a}.
$$

{\bf Proof.} It is enough to consider only the second case $ s \to \infty $
and in the first case $ s \to 0 $ under condition $ b < \infty. $ \par
{\bf A. Upper bound. } We have  for the values $ s \in (1, \infty):$
$$
\phi(s) \stackrel{def}{=} \phi(G(\psi), s) = \sup_{ p \in (a,b) }
\frac{s^{1/p} }{\psi(p) } \le s^{1/a},
$$

therefore

$$
\overline{\lim}_{s \to 0+} \frac{ \log \phi (G(\psi),s)}{ \log s } \le \frac{1}{a}.
$$

{\bf B. Low bound.} Let $ \epsilon \in (0, (b-a) ) $ be arbitrary number.
Tacking into account a simple estimation

$$
\phi(s) \ge  \frac{ s^{1/(a + \epsilon) } } {\psi(a + \epsilon)},
$$

we conclude

 $$
 \underline{\lim}_{s \to 0+} \frac{ \log \phi (G(\psi),s)}{ \log s } \ge \frac{1}{a +  \epsilon }.
$$

 This completes the proof of Lemma 3. \par

 Let now $ \vec{s} = s = (s(1), s(2), \ldots, s(k)) $ be an arbitrary {\it
positive:} $ \forall j \ \Rightarrow s(j) > 0 $ vector. We define for
arbitrary measurable function $ f: X \to R $ the following {\it family }
of linear operators (dilation operators):

$$
\sigma_s f(x) = f \left( \frac{x(d(1))}{s(1)}, \frac{x(d(2))}{s(2)}, \ldots,                          \frac{x(d(k))}{s(k)}\right).
$$
 Let also $ V_1 $ and $ V_2 $ be a two rearrangement invariant (r.i.) spaces
over $ (X, \Sigma, \mu) $  for which the operators $ \sigma_s $ are bounded
as operator from the space $ V_1 $ into the (other, in general case) space
$ V_2: $

$$
h(\vec{s}, V_1, V_2) = || \sigma_s ||(V_1 \to V_2) < \infty.
$$
{\bf Definition 1. } We define the {\it upper} particular Boyd's index
 $$
B^+_j(V_1,V_2) = B^+_j(V_1,V_2; s(1), s(2), \ldots, s(j-1), s(j+1), \ldots,
s(k) )
$$
as a limit (if there exists)

$$
B^+_j (V_1,V_2) \stackrel{def}{=} \lim_{s(j) \to \infty }
\frac{ \log h(\vec{s}, V_1, V_2)}{ \log s(j) }.
$$
Analogously may be defined the {\it low} particular Boyd's index:

$$
B^-_j(V_1,V_2) \stackrel{def}{=} \lim_{s(j) \to 0}
\frac{\log h(\vec{s}, V_1, V_2)}{ \log s(j) }.
$$

We will write for brevity for the values $ j = 1,2,\ldots, k $

$$
B^+(V_1,V_2) = \{ B^+_j(V_1,V_2) \}, \ B^-(V_1,V_2) = \{ B^-_j(V_1,V_2) \}.
$$
 These definitions may be generalized on the case when $ X = R^d $ with at
the same weight $ W(\cdot) $ and on the case when $ X = (-\pi, \pi)^d $
relatively the standard weight $ W = 1/( 2 \pi)^{d} = const. $ \par
 Notice that in the case $ d = 1, \ W(x) = 1 $ and $ V_1 = V_2 $ we obtain
the classical definition of Boyd's indices (\cite{BS1}, chapter 3, section 5,
p. 149. ) \par

{\bf  The main goal of this paper is to compute the upper and low Boyd's
(and Shimogaki's, see after) indices for pair's of Bilateral Grand Lebesgue  spaces. }  \par

 The paper is organized as follows. In the next section we calculate the
Boyd's indices for $ G(\psi; a,b) $ spaces.\par

 Further, in the section follows we compute the so \ - \ called
Shimogaki's indices in the one \ - \ dimensional case $ d = 1. $ \par

 In the fourth section we consider some generalizations of dilation operators
(multidimensional matrix dilation operators) in BGL spaces.\par

 In the last section we  study and describe  some consequences of obtained
results. \par

 We use the symbols $ C, \ C_j, \ C(X,Y), \ C(p,q;\psi) $ etc., to denote
positive finite constants along with parameters they depend on, or at least
dependence on which is essential in our study. To distinguish
between two different constants depending on the same parameters
we will additionally enumerate them, like $C_1(X,Y)$ and
$C_2(X,Y).$ The relation $ g(\cdot) \asymp h(\cdot), \ p \in (A,B), $
where $ g = g(p), \ h = h(p), \ g,h: (A,B) \to R_+, $
denotes as usually

\begin{eqnarray*}
0<\inf_{p\in (A,B)} h(p)/g(p) \le\sup_{p \in(A,B)}h(p)/g(p)<\infty.
 \end{eqnarray*}
The symbol $ \sim $ will denote usual equivalence in the limit
sense. \par
We define as usually for the two $ k \ - $ dimensional vectors $ s = \vec{s} $ and $ d = \vec{d} $
$$
s^d = \vec{s}^{\vec{d}} = \prod_{r=1}^k s(r)^{d(r)}.
$$

\bigskip

\section{Boyd's Indices}

 In this section we give an expression for the norm of a dilations operators
for  some {\it pairs} of BGL spaces  and for
the Boyd's (and other) indices for these pair of spaces.\par

\vspace{3mm}

{\bf 2.1. Dilation Operators.}\par

\vspace{3mm}

{\bf Theorem 1.1.} \\
{\bf A.} Let $ \psi(\cdot),  \zeta(\cdot), \nu(\cdot) \in E \Psi(a,b) $
and let $ \zeta = \psi \cdot \nu. $ We assert:

$$
||\sigma_{\vec{s} }||( G(\psi) \to G(\zeta)) \le \phi \left(G(\nu),
\vec{s}^{ \vec{d} + \vec{\theta} } \right). \eqno(1.A)
$$

{\bf B.} Let $ \psi(\cdot) \in \Psi, \ \zeta(\cdot), \nu(\cdot) \in E
\Psi(a,b) $ and let $ \zeta = \psi \cdot \nu. $ We assert:

$$
||\sigma_{\vec{s} }||( G(\psi) \to G(\zeta)) = \phi \left(G(\nu),
\vec{s}^{ \vec{d} + \vec{\theta} }  \right). \eqno(1.B)
$$
{\bf Proof. A.}

Let $ \psi(\cdot),  \zeta(\cdot), \nu(\cdot) \in E \Psi(a,b) $
and let $ \zeta = \psi \cdot \nu; \ f \in G(\psi). $  We can assume without
loss of generality that $ ||f||G(\psi) = 1. $  \par

 It follows from the definition of the BGL spaces  that
$$
|f|^p_{p, \mu} \le [ \ ||f||G(\psi) \ ]^p \ \psi^p(p) = \psi^p(p).
$$

We obtain using the formula of changing variables and homogeneity of the
weight $ W: $

$$
|\sigma_{\vec{s}} f|_p^p =
\int_X  \left|f \left( \frac{x(d(1))}{s(1)}, \frac{x(d(2))}{s(2)}, \ldots,                          \frac{x(d(k))}{s(k)}\right) \right|^p \ d \ \mu(x)=
$$

$$
\int_X  \left|f \left( \frac{x(d(1))}{s(1)}, \frac{x(d(2))}{s(2)}, \ldots,                          \frac{x(d(k))}{s(k)}\right) \right|^p \ W(x) \ d x =
$$

$$
\vec{s}^{\vec{d} + \vec{\theta} } \cdot
\int_X  \left|f \left( y(d(1)), y(d(2)), \ldots, y(d(k)) \right) \right|^p \ W(y) \ d y =
$$
	
$$
\vec{s}^{\vec{d} + \vec{\theta} } \cdot |f|^p_p \le \vec{s}^{\vec{d} + \vec{\theta} } \psi^p(p).
$$
 Therefore

$$
| \sigma_{\vec{s}} f|_p \le \vec{s}^{ ( \vec{d} + \vec{\theta} )/p } \ \psi(p); \ \frac{| \sigma_{\vec{s}} f|_p }{\zeta(p)} \le \frac{ \vec{s}^{ ( \vec{d} + \vec{\theta} )/p } }{ \nu(p) }.
$$

 We obtain taking the maximum over $  p \in (a,b): $

$$
|| \sigma_{\vec{s}} f||G(\zeta) = \sup_{p \in (a,b) }
\frac{| \sigma_{\vec{s}} f|_p }{\zeta(p)} \le
\sup_{p \in (a,b)} \frac{ \vec{s}^{ ( \vec{d} + \vec{\theta} )/p } }
{ \nu(p) } = \phi \left( G(\nu), \vec{s}^{\vec{d} + \vec{\theta}} \right),
$$
Q.E.D. \par

\vspace{3mm}

{\bf Proof. B.} It is easy to see that if $ \psi(\cdot) \in \Psi = \Psi(a,b), $
i.e. if $ |f|_p = \psi(p), \ p \in (a,b), $ then during all the proof of the
first assertion {\bf A } there is always the equality. \par
 Indeed,
let $ \psi(\cdot) \in \Psi, \  \zeta(\cdot), \nu(\cdot) \in E \Psi(a,b) $
and let $ \zeta = \psi \cdot \nu; \ f \in G(\psi) $ and let $ f $ be a representation of $ \psi(\cdot): \ |f|_p = \psi(p), \ p \in (a,b). $
Hence $ ||f||G(\psi) = 1. $  \par

 It follows from the definition of the BGL spaces and from the equality
 $ |f|_p = \psi(p) \ $  that
$$
|f|^p_{p, \mu} \le [ \ ||f||G(\psi) \ ]^p \ \psi^p(p) = \psi^p(p).
$$

We obtain using the formula of changing variables and homogeneity of the
weight $ W: $

$$
|\sigma_{\vec{s}} f|_p^p =
\int_X  \left|f \left( \frac{x(d(1))}{s(1)}, \frac{x(d(2))}{s(2)}, \ldots,                          \frac{x(d(k))}{s(k)}\right) \right|^p \ d \ \mu(x)=
$$

$$
\int_X  \left|f \left( \frac{x(d(1))}{s(1)}, \frac{x(d(2))}{s(2)}, \ldots,                          \frac{x(d(k))}{s(k)}\right) \right|^p \ W(x) \ d x =
$$

$$
\vec{s}^{\vec{d} + \vec{\theta} } \cdot
\int_X  \left|f \left( y(d(1)), y(d(2)), \ldots, y(d(k)) \right) \right|^p \ W(y) \ d y =
$$
	
$$
\vec{s}^{\vec{d} + \vec{\theta} } \cdot |f|^p_p = \vec{s}^{\vec{d} + \vec{\theta} } \psi^p(p).
$$
 Therefore

$$
| \sigma_{\vec{s}} f|_p = \vec{s}^{ ( \vec{d} + \vec{\theta} )/p } \ \psi(p); \
\frac{| \sigma_{\vec{s}} f|_p }{\zeta(p)} = \frac{ \vec{s}^{ ( \vec{d} + \vec{\theta} )/p } }{ \nu(p) }.
$$

 We obtain taking the upper bound over $  p \in (a,b): $

$$
|| \sigma_{\vec{s}} f||G(\zeta) = \sup_{p \in (a,b) }
\frac{| \sigma_{\vec{s}} f|_p }{\zeta(p)} =
\sup_{p \in (a,b)} \frac{ \vec{s}^{ ( \vec{d} + \vec{\theta} )/p } }
{ \nu(p) } = \phi \left( G(\nu), \vec{s}^{\vec{d} + \vec{\theta}} \right),
$$

 This completes the proof of Theorem 1.1. \par

\vspace{3mm}

{\bf 2.2. Boyd's Indices. Main Result.} \par

\vspace{3mm}

{\bf Theorem 1.2.} Let again $ \psi(\cdot) \in \Psi(a,b), \ \zeta(\cdot), \nu(\cdot) \in E \Psi(a,b) $ and let $ \zeta = \psi \cdot \nu. $ We assert:

$$
B^+(G(\psi), G(\zeta) ) = \frac{ \vec{d} + \vec{\theta} }{a},
\eqno(1.C)
$$

$$
B^-(G(\psi), G(\zeta) ) = \frac{ \vec{d} + \vec{\theta} }{b}.
\eqno(1.D)
$$

{\bf Proof} it follows from the assertion of theorem 1.1 for the explicit
view of Boyd's indices and from the Lemma 3. Namely,

$$
B^+_j (G(\psi),G(\zeta)) = \lim_{s(j) \to \infty }
\frac{ \log h(\vec{s}, G(\psi), G(\zeta))}{ \log s(j) } =
$$

$$
\lim_{s(j) \to \infty } \frac{ \log \phi ( G (\nu), s^{d + \theta} ))}
{ \log s(j) } = \frac{ d(j) + \theta(j)}{a}.
$$

\vspace{3mm}

{\bf Corollary.} At the same assertions as in the theorems 1.1 and 1.2 are
true for the spaces $ G^o(\zeta), \ G^o(\psi), \ G^o(\nu) $ and following
for the spaces $ GA(\zeta), \ GA(\psi), \ GA(\nu); \ GB(\zeta), \ GB(\psi), \ GB(\nu). $ \par

 Namely, \par
{\bf AA.} Let $ \psi(\cdot),  \zeta(\cdot), \nu(\cdot) \in E \Psi(a,b) $
and let $ \zeta = \psi \cdot \nu. $ We assert:

$$
||\sigma_{\vec{s} }||( G^o(\psi) \to G^o(\zeta)) \le \phi \left(G^o(\nu),
\vec{s}^{ \vec{d} + \vec{\theta} } \right); \eqno(1.AA)
$$

{\bf BB.} Let $ \psi(\cdot) \in \Psi(a,b), \ \zeta(\cdot), \nu(\cdot) \in E
\Psi(a,b) $ and let $ \zeta = \psi \cdot \nu. $ We assert:

$$
||\sigma_{\vec{s} }||( G^o(\psi) \to G^o(\zeta)) =
\phi \left(G^o(\nu),\vec{s}^{ \vec{d} + \vec{\theta} }  \right) =
$$

$$
\phi \left(G(\nu),\vec{s}^{ \vec{d} + \vec{\theta} }  \right); \eqno(1.BB)
$$

$$
B^+(G^o(\psi), G^o(\zeta) ) = \frac{ \vec{d} + \vec{\theta} }{a};
\eqno(1.CC)
$$

$$
B^-(G^o(\psi), G^o(\zeta) ) = \frac{ \vec{d} + \vec{\theta} }{b}.
\eqno(1. DD)
$$
{\bf Proof.} It is enough to prove that if $ d = 1 $ and $ \psi \in \Psi, \ \zeta \in E \Psi, \nu \in E \Psi, $ then

$$
|| \sigma_{ \vec{s} } || \left( G^o(\psi) \ \to \ G^o(\zeta) \right) =
\phi \left(G(\nu),\vec{s}^{ \vec{d} + \vec{\theta} }  \right).
$$
 Let $ f: X \to R $ be the measurable {\it non \ - \ negative }
 function  for which $ |f|_p = \psi(p),
\ p \in (a,b). $ The existence of this function it follows from the condition
$ \psi \in \Psi. $  We have from the definition of the norm of a function
belonging to the space $ G(\psi) $ or $ G^o(\psi) $ that $ ||f||G(\psi) = 1. $
\par

 Since the  measure $ \mu $ is sigma \ - \ finite, there exists a decreasing:
$ A(1) \subset A(2) \subset A(3) \ldots $
sequence of measurable sets $ A(n), n = 1,2,\ldots, \ A(n) \in \Sigma $ with finite measures $ \mu(A(n)) < \infty $ covering the space $ X: $

$$
\cup_{n=1}^{\infty} A(n) = X.
$$
 Let us consider the sequence of a non \ - \ negative functions
$$
f_n = f_n(x) = f(x) \cdot I(A(n),x) \cdot I( x: |f(x)| \le n) ).
$$

 We observe: $ f_n \uparrow f, \ f_n \in G^o(\psi) $ as long as $ G^o(\psi)  = GB(\psi), \ ||f_n||G^o(\psi) = ||f_n||G(\psi) \uparrow ||f||G(\psi),
\ \forall p \in (a,b) \ \Rightarrow \ |f_n|_p \uparrow |f|_p. $  \par
 The inequality (1.AA) is obvious; let us prove the inverse inequality. Since
for all sufficient great values $ n $

$$
|| \sigma_{ \vec{s} } || \left( G^o(\psi) \ \to \ G^o(\zeta) \right) \ge
\frac{ ||\sigma_s f_n||G^o(\zeta)}{||f_n||G^o(\psi) } \ge
\frac{ ||\sigma_s f_n||G(\zeta)}{||f||G(\psi) }.
$$

We get using Fatou property of $ G(\psi) $ spaces:

$$
|| \sigma_{ \vec{s} } || \left( G^o(\psi) \ \to \ G^o(\zeta) \right) \ge
\sup_n \frac{ ||\sigma_s f_n||G(\zeta)}{||f||G(\psi) } =
$$

$$
\sup_n \frac{ \sup_{p \in (a,b) } |\sigma_s f_n|_p /\zeta(p)}{||f||G(\psi) } =
\sup_{p \in (a,b) } \frac{ \sup_n  |\sigma_s f_n|_p /\zeta(p)}{||f||G(\psi) } =
$$

$$
\sup_{p \in (a,b) } \frac{ |\sigma_s f|_p /\zeta(p)}{||f||G(\psi) } =
|| \sigma_{ \vec{s} } || \left( G(\psi) \ \to \ G(\zeta) \right)	=
\phi \left( G(\nu), \vec{s}^{\vec{d} + \vec{\theta}} \right).
$$

 For example, if $ \zeta(p) = \psi(p), \ p \in (a,b); \ d = 1, \ W(x) = 1; $
then (cf. \cite{Holon})

$$
B^+(G(\psi), G(\psi)) = B^+(G^o(\psi), G^o(\psi)) = 1/a,
$$

$$
 B^-(G(\psi), G(\psi)) = B^+(GA(\psi), GB(\psi)) = 1/b.
$$

\bigskip

\section{Shimogaki's indices }

\vspace{3mm}

 We consider in this section only the one \ - \ dimensional case $ d = 1 $. \par

 Let $ G $ be again r.i. space over $ (X, \Sigma, \mu ) $ with correspondent
fundamental function $ \phi(G, \delta) = \phi(\delta). $  Let us denote

$$
M_G(t) = \sup_{s > 0} \frac{\phi(s t) }{\phi(s)},
$$

$$
\beta^-(G) = \sup_{t \in (0,1) } \frac{ \log M_G(t)}{\log t} =
 \lim_{t \to 0+ } \frac{ \log M_G(t)}{\log t},
$$

$$
\beta^+(G) = \inf_{t \in (1, \infty) } \frac{ \log M_G(t)}{\log t} =
 \lim_{t \to \infty } \frac{ \log M_G(t)}{\log t}.
$$
 The numbers $ \beta^-(G) $ and $ \beta^+(G) $ are called
Shimogaki indices. It is known (see \cite{ BS1}, p. 171 \ - \ 178) that

$$
0 \le B^-(G,G) \le \beta^-(G) \le \beta^+(G) \le  B^+(G,G).
$$

\vspace{2mm}

{\bf Theorem 2. } {\it Let } $ G = G(\psi; a,b), \ \psi \in \Psi; $ {\it
then }

$$
1/b =  \beta^-(G(\psi; a,b)) < \beta^+(G(\psi; a,b)) = 1/a.
$$

\vspace{2mm}

 {\bf Proof.}\\
{\bf A. Upper bound.}  Assume for certainty $ t \to \infty, \ t > 1. $

$$
\phi(st) = \sup_{p \in (a,b)} \frac{ s^{1/p} \ t^{1/p} }{\psi(p) } \le
 t^{1/a} \  \sup_{p \in (a,b)} \frac{ s^{1/p} }{\psi(p) } = t^{1/a} \phi(s),
$$

hence

$$
\overline{\lim}_{t \to \infty } \frac{ \log M_G(t)}{\log t} \le 1/a.
$$

{\bf B. Low bound.}  Let $ \epsilon = const \in (0, (b-a)). $

$$
M_G(t) \ge \frac{\phi((a + \epsilon)t ) }{ \phi(a + \epsilon) } \ge
C(a,\epsilon) \ \phi((a + \epsilon)t).
$$

 We find using Lemma 3:

$$
\underline{\lim}_{t \to \infty } \frac{ \log M_G(t) }{ \log t } \ge
1/(a + \epsilon).
$$
 Since the number $ \epsilon $ is arbitrary, we obtained the proof
what is desired. \par

\vspace{2mm}

{\bf Corollary.} The condition of homogeneity of the weight $ W = W(\vec{x}) $
may be  weakened as follows. \par
 Let us denote

$$
K^+_{\infty} = \sup_{ \{ \min \ s(j) \ge 1 \}} \frac{ W(s y)}
{ s^{\theta} \ W(y)},
$$

$$
K^-_{\infty} = \inf_{ \{ \min \ s(j) \ge 1 \}} \frac{ W(s y)}
{ s^{\theta} \ W(y)},
$$

$$
K^+_{0} = \sup_{ \{ \max \ s(j) \le 1 \}} \frac{ W(s y)}
{ s^{\theta} \ W(y)},
$$

$$
K^-_{0} = \inf_{ \{ \max \ s(j) \le 1 \}} \frac{ W(s y)}
{ s^{\theta} \ W(y)}.
$$
 We assert: \par
{\sc I.} If $  K^+_{\infty} < \infty, \ \psi, \zeta,  \nu  \in E\Psi(a,b),
\ \min(s(j)) > 1, $ then

$$
||\sigma_{\vec{s} }||( G(\psi) \to G(\zeta)) \le \max
\left[ (K^+_{\infty})^{1/a}, \ (K^+_{\infty})^{1/b} \right] \cdot
\phi \left(G(\nu),\vec{s}^{ \vec{d} + \vec{\theta} } \right);
$$

{\sc II.} If $  K^-_{\infty} < \infty, \ \psi \in \Psi(a,b), \zeta,  \nu
\in E\Psi(a,b), \  \min(s(j)) > 1, $ then

$$
||\sigma_{\vec{s} }||( G(\psi) \to G(\zeta)) \ge \min
\left[ (K^-_{\infty})^{1/a}, \ (K^-_{\infty})^{1/b} \right] \cdot
\phi \left(G(\nu),\vec{s}^{ \vec{d} + \vec{\theta} } \right);
$$

{\sc III.} If $  K^+_{0} < \infty, \ \psi, \zeta,  \nu
\in E\Psi(a,b), \  \max(s(j)) > 1, $ then

$$
||\sigma_{\vec{s} }||( G(\psi) \to G(\zeta)) \le \max
\left[ (K^-_{0})^{1/a}, \ (K^-_{0})^{1/b} \right] \cdot
\phi \left(G(\nu),\vec{s}^{ \vec{d} + \vec{\theta} } \right);
$$

{\sc IV.} If $  K^-_{0} > 0, \ \psi \in \Psi(a,b), \zeta,  \nu
\in E\Psi(a,b), \  \max(s(j)) > 1, $ then

$$
||\sigma_{\vec{s} }||( G(\psi) \to G(\zeta)) \ge \min
\left[ (K^-_{0})^{1/a}, \ (K^-_{0})^{1/b} \right] \cdot
\phi \left(G(\nu),\vec{s}^{ \vec{d} + \vec{\theta} } \right);
$$
 As a consequence: if all the conditions { \sc I \ - \ IV } are satisfied, then

$$
B^+(G(\psi), G(\zeta) ) = \frac{ \vec{d} + \vec{\theta} }{a},
$$

$$
B^-(G(\psi), G(\zeta) ) = \frac{ \vec{d} + \vec{\theta} }{b}.
$$
 The conclusion of theorem 2 also holds. \par
\bigskip

\section{Matrix Dilations}

\vspace{3mm}

 In this section we describe briefly some multidimensional matrix dilation
operators in BGL spaces and compute its norms.\par
 Let $ X = R^d $ be usually Euclidean space with Lebesgue measure $ m:
m(dx) = dx. $ Let $ A $ be a square: $ A = d \times d $   non \ - \ degenerate constant matrix. We define
alike  in \cite{Olphert} the matrix dilation operator $ D_A $ as follows:
for arbitrary measurable function $ f: X \to R $

$$
(D_A \ f)(x) = f \left(A^{-1} \ x \right).
$$

 As

$$
| D_A \ f|^p_p = \int_X |f \left( A^{-1} \ x \right)|^p \ dx =
$$
$$
| det(A)| \ \int_X  |f(y)|^p \ dy  = | det(A)| \ |f|_p^p,
$$
we conclude in the case $ \psi(\cdot), \zeta(\cdot), \nu(\cdot) \in
E \Psi( a,b), \ \zeta(p) = \psi(p) \cdot \nu(p):$

$$
\frac{|D_A \ f|_p}{\zeta(p)} \le ||f||G(\psi) \ \frac{|det(A)|^{1/p}}
{\nu(p)};
$$

$$
||D_A \ f||G(\zeta) \le ||f||G(\psi) \cdot \sup_{p \in (a,b)}
\frac{|det(A)|^{1/p}}{ \nu(p)} =
$$

$$
||f||G(\psi) \cdot \phi(G(\nu), |det(A)|).
$$

Therefore,

$$
||D_A||(G(\psi) \to G(\zeta) ) \le  \phi(G(\nu), |det(A)|)
$$

with equality  in the case, e.g. if
 $ \psi(\cdot) \in \Psi(a,b), \zeta(\cdot), \nu(\cdot) \in
E \Psi(a,b), \ \zeta(p) = \psi(p) \cdot \nu(p).$ \par

We have as a consequence in the last case:

$$
\lim_{det(A) \to \infty} \frac{\log ||D_A||(G(\psi) \to G(\zeta) )}
{\log |det(A)|} = \frac{1}{a},
$$

$$
\lim_{det(A) \to 0} \frac{\log ||D_A||(G(\psi) \to G(\zeta) )}
{\log |det(A)|} = \frac{1}{b}.
$$
 We consider now some slight generalization. Let  $ |||x||| = ||| \vec{x}||| $
be some non \ - \ degenerate norm in the space $ R^d $ and let $ |||A||| $ be
the correspondence norm of the matrix $ A:$

$$
|||A||| = \sup \{ |||A \ x|||, \ x: |||x||| = 1 \}.
$$

 Let $ \mu_{\sigma} $ be a following weight measure:

$$
\mu_{\sigma}(V) = \int_V |||x|||^{\sigma} \ \ dx, \ \sigma  = const.
$$
 We assert in "general" case $ \psi(\cdot), \zeta(\cdot), \nu(\cdot) \in
E \Psi( a,b), \ \zeta(p) = \psi(p) \cdot \nu(p):$

$$
||D_A||( G(\psi, \mu_{\sigma}) \to G(\zeta, \mu_{\sigma}) ) \le
\phi(G(\nu, \mu_{\sigma}), |det(A)| \cdot |||A|||^{\sigma} ).
$$

We get as an particular case, i.e.
in the case if in addition the norm $ ||| \cdot||| $ is usually Euclidean
norm and $ A = s \times U, $ where $ s = const > 0 $ and $ U $ be the
unitarian  operator:

$$
||D_A||(G(\psi, \mu_{\sigma}) \to G(\zeta, \mu{\sigma}) ) \le
\phi(G(\nu, \mu_{\sigma}), s^{d + \sigma}).
$$

Therefore, we conclude in the considered case under additional assumptions
$ \psi(\cdot) \in \Psi(a,b), \zeta(\cdot), \nu(\cdot) \in
E \Psi( a,b), \ \zeta(p) = \psi(p) \cdot \nu(p):$

$$
\lim_{s \to \infty}  \frac{\log ||D_A||(G(\psi) \to G(\zeta) )}
{\log s} = \frac{d + \sigma}{a},
$$

$$
\lim_{s \to 0}  \frac{\log ||D_A||(G(\psi) \to G(\zeta) )}
{\log s} = \frac{d + \sigma}{b}.
$$

 The case when the operator $ A $ is the Croneker product of some non \ - \
degenerate  linear operators may be considered analogously.  See for the
definitions and preliminary results, e.g., \cite{Olphert}. \par

\bigskip

\section{Concluding Remarks }

\vspace{3mm}

 We will show here some applications of obtained results. \par
 In this section we consider  only the one \ - \  dimensional case $ d = 1, $
i.e.  the cases $ X = R, \ X = R^+ \ $ or  $ X = (0, 2 \pi); $ with Lebesgue measure:
$ W(\vec{x}) = 1, $ i.e. $ \vec{\theta} = 0. $ \par

\vspace{2mm}

{\bf 1.} {\it Conjugate spaces.} \par
 The {\it associate } spaces $ (G(\psi))^' $  to the BGLS $ (G(\psi)) $
are described in \cite{Holon}. Using
the Corollary 4.2 from \cite{BS1}, chapter 1, section 4, we compute the
{\it conjugate } spaces $ (GA(\psi))^* $ to the $ GA(\psi) $ spaces:

$$
(GA(\psi))^* = (GB(\psi))^* = (G^o(\psi))^* =  (G(\psi))^/.
$$

\vspace{2mm}

{\bf 2.} {\it Boyd indices for associate spaces. } \par

\vspace{2mm}

It follows from \cite{BS1}, chapter 3, section 5, proposition 5.13 that
if $ \psi \in \Psi(a,b), \ 1 \le a < b \le \infty, $ then
$$
B^+( (G(\psi))^/,  (G(\psi))^/ ) = 1 \ - \ 1/b,
$$

$$
B^-( (G(\psi))^/,  (G(\psi))^/ ) = 1 \ - \ 1/a.
$$

\vspace{2mm}

{\bf 3.} {\it Boundedness of some singular operators.}\par

\vspace{2mm}

{\sc I. } Let $ X = R^1_+, \  \psi \in \Psi(a,b), \ 1 \le a < b \le \infty, $
and consider two linear operators of a Hardy \ - \ Littlewood type:

$$
(P_{\alpha})f(t) = t^{\ - \ \alpha} \int_0^t s^{\alpha \ - \ 1} \ f(s) \ ds,
$$

$$
(Q_{\beta}) f(t) = t^{ \ - \ \beta} \int_t^{\infty} s^{ \beta \ - \ 1}
\ f(s) \ ds,
$$

$ s,t \in (0,\infty), \ \alpha, \beta = const \in (0,1). $ \par
We conclude using the theorem 5.15 from \cite{BS1}, chapter 3, section 5 that
the operator $ P_{\alpha} $ is bounded in the space $ G(\psi; a,b) $ iff
$ \alpha > 1/a; $ the operator $ Q_{\beta} $ is bounded in the space
$ G(\psi; a,b) $ iff $ \beta < 1/b. $ \par

{\sc II. } Let $ X = R^d_+, \ W(x) = 1,  \psi \in \Psi(a,b), \ 1 \le a < b \le \infty,  $ and consider  the (quasi \ - \ linear) Hardy  \ - \ Littlewood
maximal operator $ M: $

$$
(Mf)(x) = \sup_Q  \left[ \int_Q |f(y)| \ dy \ /|Q| \right],
$$
where the supremum  extends over all non \ - \ degenerate cubes $ Q $
containing $ x $ (cubes will be assumed to have  their sides parallel to the coordinate axes), $ |Q| $ denotes the $ d \ - $ dimensional volume of $ Q. $ \par

 We conclude using the theorem 5.17, belonging to G.G.Lorentz and T.Shimogaki,
 from \cite{BS1}, chapter 3, section 5 that
the operator $ M $ is bounded in the space $ G(\psi; a,b) $ iff
$ a > 1. $ \par

{\sc III. } Let $ X = R^1, \ W(x) = 1,  \psi \in \Psi(a,b), \ 1 \le a < b \le \infty,  $ and consider  the Hilbert transform $ H. $ \par
 We conclude using the theorem 5.18, belonging to  D.W.Boyd,
from \cite{BS1}, chapter 3, section 5 that
the operator $ H: $

$$
(Hf)(x) = \pi^{-1} \ \lim_{\epsilon \to 0+} \int_{ \{y: |x \ - \ y| >
\epsilon \} } f(y) \ \frac{ d y }{ x - y }
$$

is bounded in the space $ G(\psi; a,b) $ iff $ a > 1, \ b < \infty. $ \par

{\bf 4.} {\it Norm convergence of Fourier series.} \par

Let here $ X = (0, 2 \pi), \ W(x) = 1,  \psi \in \Psi(a,b), \ 1 \le a < b \le
\infty. $ We consider the usual Fourier series for {\it arbitrary } function
$ f: \ X \to R $  and such that $ f \in  G^o(\psi) ( = GA(\psi) = GB(\psi)). $
\par
 We obtain using the corollary 6.11 from  \cite{BS1}, chapter 3, section 6
that the Fourier series for arbitrary function $ f  \in G^o(\psi) $ converge
in the norm $ G(\psi; a,b) $ if and only if

$$
a > 1, \ b < \infty.
$$

 Note in addition to this section that the other cases, e.g. if $ a = 1 $
or/and $ b = \infty $ are complete investigated in \cite{Holon}.  In this case
the considered singular operators: Hilbert, Hardy \ - \ Littlewood, Fourier etc.
are bounded as an operators from one BGL space into {\it another } space. \par

\vspace{3mm}

{\sc Department of Mathematics, Bar \ - \ Ilan University, 59200, Ramat
\ - \ Gan,  ISRAEL. }\\

\vspace{3mm}

 E \ - \ mail address: \ eugeny@soniclynx.net \\

\vspace{3mm}

 E \ - \ mail address: \ sirota@zahav.net.il \\

\end{document}